# INTERPOLATION BY SUMS OF SERIES OF EXPONENTIALS AND GLOBAL CAUCHY PROBLEM FOR CONVOLUTION OPERATORS

S. G. Merzlyakov[1,*] and S. V. Popenov[2,*]

INTERPOLATION PROBLEM

Let $M = \{\mu_k\} \subset \mathbb{C}$ be a countable nodes set which is discrete in $\mathbb{C}$; $m_k \in \mathbb{N}$ be the multiplicities of nodes $\mu_k$; and $\Lambda \subset \mathbb{C}$ be some unbounded exponents set. Let us define $\Sigma(\Lambda, \mathbb{C}) = \{U: U(z) = \sum_{n=1}^{\infty} c_n e^{\lambda_n}, \lambda_n \in \Lambda\}$ where the series is absolutely convergent for $z \in \mathbb{C}$. It is known that such the series converges uniformly on all compact subsets of $\mathbb{C}$, so every $U \in \Sigma(\Lambda, \mathbb{C})$ is an entire function; $\Sigma(\Lambda, \mathbb{C})$ is a linear subspace of $H(\mathbb{C})$ that in general it is not the closed one.

The aim is to describe classes of sets $M$ and $\Lambda$ that give solubility of the problem: *for arbitrary interpolation data $b_k^j \in \mathbb{C}, k \in \mathbb{N}, j = 1, \ldots, m_k - 1$, there exists such a sum $U \in \Sigma(\Lambda, \mathbb{C})$ that $U^{(j)}(\mu_k) = b_k^j$.*

Counterexample. For $\Lambda = \{2\pi n, n \in \mathbb{N}\}$ all entire functions of the form $U(z) = \sum_{n=1}^{\infty} c_n e^{2\pi n z}$ are solutions of the difference equation $U(z) = U(z + i)$. Thus for any two nodes $\mu_1 = \mu$, $\mu_2 = \mu + i$, $\mu \in \mathbb{C}$, the simple interpolation $U(\mu_k) = b_k$ fails for interpolation data $b_1 \neq b_2$.

Representations for $H(\mathbb{C})$. The classical result on interpolation in $H(\mathbb{C})$ implies the following. Let $\psi_M$ be a function in $H(\mathbb{C})$ whose zero set is $M = \{\mu_k\}$ with multiplicities $m_k$ accounted. Let $(\psi_M) = \{h \in H(\mathbb{C}): h = r \cdot \psi_M, r \in H(\mathbb{C})\}$ be the ideal in $H(\mathbb{C})$ generated by $\psi_M$. For a given nodes set $M$ solubility of the interpolation problem is equivalent to existence of representation $H(\mathbb{C}) = \Sigma(\Lambda, \mathbb{C}) + (\psi_M)$.

Definition. Denote by $\mathbb{S} = \{s \in \mathbb{C}: |s| = 1\}$. Let $X \subset \mathbb{C}$ be an unbounded set. The set $P(X) \subset \mathbb{S}$ of limit directions at infinity consists of all such $s \in S$ that $s = \lim_{k \to \infty} x_k /$

---
[1] Institute of Mathematics, Ufa Federal Research Centre, Russian Academy of Sciences, Ufa
[*] msg@mail.ru
[2] Institute of Mathematics, Ufa Federal Research Centre, Russian Academy of Sciences, Ufa
[*] spopenov@gmail.com

$|x_k|$ for some sequence $\{x_k\} \subset X$, $x_k \to \infty$.

Let us select in arbitrary way a sparse sequence $\{\lambda_n\} \subset \Lambda$ such that $|\lambda_{n+1}| > 2|\lambda_n|$ and $P(\{\lambda_n\}) = P(\Lambda)$. It has zero density and entire function $G_1(z) = \prod_{n=1}^{\infty}(1 - \frac{z}{\lambda_n})$ has minimal type for order of growth 1.

Let $M_{G_1}$ be the convolution operator in $H(\mathbb{C})$ generated by $G_1$. The results of A. F. Leont`ev [1] (as well as O. A. and A. S. Krivosheev [2]-[4]) imply the Fundamental Principle for $M_{G_1}$ that is $Ker M_{G_1} = \Sigma(\{\lambda_n\}, \mathbb{C})$. Thus one should study the representations $H(\mathbb{C}) = Ker M_{G_1} + (\psi_M)$ because $Ker M_{G_1} = \Sigma(\{\lambda_n\}, \mathbb{C}) \subset \Sigma(\Lambda, \mathbb{C})$.

GLOBAL CAUCHY PROBLEM FOR CONVOLUTION OPERATORS

When $D = \mathbb{C}$ and $M$ is a discrete subset of some angle V.V. Napalkov [5]-[7] is the first who obtains some sufficient conditions for existence of the representations

$$H(\mathbb{C}) = Ker M_G + (\psi_M), \qquad (1)$$

where $M_G$ is the convolution operator generated by an entire function $G$ of exponential type. Every element of $Ker M_G$ is a limit of quazi-polynomials whose exponents are from the zero set $\Lambda = Z_G$ of function $G$ but in general it is not represented by any series of exponentials. As specified above, in representations (1) one may replace the kernel $Ker M_G$ on $Ker M_{G_1} = \Sigma(\{\lambda_n\}, \mathbb{C})$, where $\{\lambda_n\} \subset Z_G$ where $Z_G$ is a sparse sequence.

We introduce the special class of sets $M$ that includes the ones considered before. Then the conventional duality is used which is given by the Laplace transformation of the strong dual space and the method is found to treat natural dual problems in the space of entire functions of exponential type. That enable us to obtain conditions of solubility of interpolation and as the result of the global Cauchy problem with data on the special $M$ in the form of series of exponentials, at the same time we may treat holomorphic function in a domain. This strengthens part of results of V.V. Napalkov and his disciples (see e.g. [5]-[8] as well as our results in [9]-[11]).

It turns out that for some $\Lambda$ the obtained conditions are the necessary ones for more general interpolation by functions of the form $U(z) = \int_\Lambda e^{\lambda z} d\sigma$ where $d\sigma$ is arbitrary Radon measure on $\Lambda$ and in addition for arbitrary $M$. Meanwhile the research of necessary conditions for the Cauchy problem (1) is an open complicated problem.

# INTERPOLATION DEFECT

We use hereafter the notation $s_\tau = e^{i\tau}$ and $s_\omega = e^{i\omega}$, $\omega, \tau \in [0, 2\pi)$, to refer to elements of $P(M)$ and $P(\Lambda)$ respectively.

Conditions. For every $s_\tau \in P(M)$ one of two following conditions is fulfilled:

(A) A direction $s_\omega \in P(\Lambda)$ exists such that $\mathrm{Re}(s_\omega s_\tau) > 0$ and for a number $\varkappa = \varkappa(s_\omega) \in \mathbb{R}$ and for all $c > \varkappa$ every line $l(s_\omega, c) = \{z \in \mathbb{C} : \mathrm{Re}(s_\omega z) = c\}$ contains at most one node $\mu_k \in M$.

(B) A direction $s_\omega \in P(\Lambda)$ exists such that $\mathrm{Re}(s_\omega s_\tau) = 0$ and for a number $\varkappa = \varkappa(s_\omega) \in \mathbb{R}$ and for all $c > \varkappa$ every line $l(s_\omega, c)$ contains at most one node $\mu_k \in M$; for every subsequence $\{\mu_{k_l}\} \subset M$ such that $|\mu_{k_l}| \to \infty$ and $s_\tau = \lim\limits_{k_l \to \infty} \mu_{k_l} / |\mu_{k_l}|$ there exists limit $\lim\limits_{k_l \to \infty} \mathrm{Re}(s_\omega \mu_{k_l}) = +\infty$.

Denote by $P_M(\Lambda) \subset P(\Lambda)$ the set of all such limit directions $s_\omega \in P(\Lambda)$ that are coupled with $s_\tau \in P(M)$ in (A) and (B). The conditions implies that the set $D_M = \bigcap_{s_\omega \in P_M(\Lambda)} \{\mu_k \in M : \mathrm{Re}(s_\omega \mu_k) \leq \varkappa(s_\omega)\}$ is finite or empty.

Consider all quasi-polynomials $P(z)$ generated by monomials $z^j \exp(\mu_k z)$ where $\mu_k \in D_M$, $j = 0, \ldots, m_k - 1$. Let $d_M \in \mathbb{Z}_+$ be the dimension of subspace of all $P$ that vanish on $\{\lambda_n\}$. In some cases $d_M = 0$; if $D_M \neq \emptyset$ the condition means that $\{\lambda_n\}$ is the uniqueness set for all such quasi-polynomials $P$.

Теорема 1. *Assume that conditions (A) and (B) are fulfilled then for arbitrary sparse sequence $\{\lambda_n\} \subset \Lambda$, i.e. $|\lambda_{n+1}| > 2|\lambda_n|$, $P(\{\lambda_n\}) = P_M(\Lambda)$,*

*1) subspace $\mathrm{Ker}\, M_{G_1} + (\psi_M)$, is closed; if $d_M > 0$ it has finite codimension $d_M$ in $H(\mathbb{C})$; moreover there exists such $\xi_1, \ldots, \xi_{d_M} \in \mathbb{C}$ that $H(\mathbb{C}) = \mathrm{Ker}\, M_{G_1} + (\psi_M) + \mathrm{span}\, \{\exp(\xi_l z)\}_{l=1}^{d_M}$;*

*2) if $d_M = 0$ one has $H(\mathbb{C}) = \mathrm{Ker}\, M_{G_1} + (\psi_M)$;*

*3) subspace $\mathrm{Ker}\, M_{G_1} \cap (\psi_M)$ has infinite dimension.*

Claim 2) gives solubility of considered interpolation problem in the domain (as in corollary 1 below) and of the Cauchy problem (1) for general convolution operator $M_G$ (recall that in the last case $\mathrm{Ker}\, M_{G_1} = \Sigma(\{\lambda_n\}, \mathbb{C}) \subset \mathrm{Ker}\, M_G, \{\lambda_n\} \subset Z_G$) when $D_M = \emptyset$. If $D_M \neq \emptyset$ one may give readily verified tests for $d_M = 0$. E.g. $d_M = 0$ provided there exists a direction $s_\omega$ such that every line $l(s_\omega, c)$ contains at most one node $\mu_k \in M$. Claim 3) indicates on the structure of non-uniqueness for interpolation.

Using well-known results (see e.g. [12]) on asymptotic distribution and localization of zeroes of quasi-polynomials we can give more general tests both for closedness of the sum and for $d_M = 0$.

## MORE GENERAL APPROACH TO INTERPOLATION

Let us show that one may employ more general objects in study of interpolation problem. As an illustration, suppose that $\Lambda$ has a unique limit direction $s_\omega$ at infinity.

Let us modify the class of considered nodes sets: in what follows sets $M \subset \mathbb{C}$ may have finite limit points but every $M$ does not contain them. Additionally let us upgrade the definition of series of exponentials employed for interpolation on the base of the known analogue of the Abel theorem for power series [13].

**The Abel theorem.** *If exponents of the series of exponentials have the unique limit direction $s_\omega$ and it absolutely converges in a neighborhood of a point $z = \mu_k$ the series absolutely converges in the half-plane $\Pi(s_\omega, \mu_k) = \{z \in \mathbb{C}: Re(s_\omega z) < Re(s_\omega \mu_k)\}$.*

The theorem enable us to assume that the series of exponentials employed for interpolation are the ones absolutely converging merely in neighborhoods of all $\mu_k \in M$ hence their sums are locally analytic functions on $M$. But any such a series is absolutely convergent in domain $\Pi = \Pi_{\omega,M} = \{z \in \mathbb{C}: Re(s_\omega z) < d(s_\omega, M) = \sup_{\mu_k \in M}(s_\omega \mu_k)\}$ and its sum is holomorphic in $\Pi$.

If $d(s_\omega, M) = +\infty$ the domain $\Pi$ is the plane $\mathbb{C}$ and as above the sums of series from $\Sigma(\Lambda, \mathbb{C})$ are employed for interpolation; the sums are the entire functions. If $d(s_\omega, M) < +\infty$ the domain $\Pi$ is the half-plane $\{z \in \mathbb{C}: Re(s_\omega z) < d(s_\omega, M)\}$ and in that case according to the new definition the sums of series from $\Sigma(\Lambda, \Pi)$ that are absolutely convergent in $\Pi$ are employed for interpolation. The sums are holomorphic functions in the half-plane $\Pi$ so one should study representations $H(\Pi) = \Sigma(\Lambda, \Pi) + (\psi_M)_\Pi$ where $\psi_M \in H(\Pi), Z_{\psi_M} = M$.

## CRITERIA FOR SPECIAL CLASS OF NODES SETS

In what follows, along with the specified changes in the set up of the problem, nodes sets $M$ satisfy more strong condition then in (A) and (B) above: for $c \in \mathbb{R}$ every straight line $l(s_\omega, c) = \{z \in \mathbb{C}: Re(s_\omega z) = c\}$ contains at most one node $\mu_k \in M$.

Теорема 2. *Assume that $P(\Lambda) = \{s_\omega\}$ and let $\{\lambda_n\} \subset \Lambda$ be any sparse sequence.*
1. *Let $d(s_\omega, M) = +\infty$ i.e. $\Pi = \mathbb{C}$. The problem of interpolation by sums of series from $\Sigma(\{\lambda_n\}, \mathbb{C})$ is soluble in $H(\mathbb{C})$ for all special sets $M$ described above if and only if the limit exists, $\lim_{k \to \infty} Re(s_\omega \mu_k) = +\infty$.*
2. *Let $d(s_\omega, M) < \infty$ i.e. $\Pi$ is the half-plane. The problem of interpolation by sums of series from $\Sigma(\{\lambda_n\}, \Pi)$ is soluble in $H(\Pi)$ for all special sets $M$ described above if and only if $\text{Re}(s_\omega \mu_k) < d(s_\omega, M)$ for all $\mu_k \in M$ and $\lim_{k_i \to \infty} \text{Re}(s_\omega \mu_{k_i}) = d(s_\omega, M)$ exists.*

Corollary. *Let us assume that $P(\Lambda) = \{s_\omega\}$; let $D$ be an arbitrary domain and $M \subset D$ be a special discrete set in $D$ described above. The problem of interpolation by sums of series absolutely converging in neighborhoods of all $\mu_k \in M$ is soluble in $H(D)$ if and only if all finite limit points of $M$ belong to $\partial D \cap \Pi$ and the conditions of claims 1 and 2 of theorem 1 are satisfied. The interpolation is realized by sums of series from $\Sigma(\{\lambda_n\}, \Pi) \subset H(D)$ where $\{\lambda_n\} \subset \Lambda$ is any sparse sequence.*

Crucially, the necessity in theorem 2 has been proved for arbitrary nodes sets and, as mentioned above, for general problem of interpolation by the Laplace transforms of the Radon measures on $\Lambda$. For brevity of notation we consider herein the special sets $M$ and the sets $\Lambda$ with unique limit direction. Our methods enable us to study more general sets $M$ and $\Lambda$ as well.